\documentclass[12pt]{article}
\usepackage{amssymb}
\usepackage{amsmath}
\usepackage{latexsym}
\usepackage{amsthm}

\makeatletter
\def\thmhead@plain#1#2#3{%
  \thmname{#1}\thmnumber{\@ifnotempty{#1}{ }#2}%
  \thmnote{ \the\thm@notefont(#3)}}
\let\thmhead\thmhead@plain
\def\swappedhead#1#2#3{%
  \thmnumber{#2}\thmname{\@ifnotempty{#2}{. }#1}%
  \thmnote{ \the\thm@notefont(#3)}}
\makeatother

\theoremstyle{definition} 

 \newtheorem{definition}{Definition}[section]
 \newtheorem{remark}[definition]{Remark}

\theoremstyle{plain}      

 \newtheorem{proposition}[definition]{Proposition}
 \newtheorem{theorem}[definition]{Theorem}
 \newtheorem{corollary}[definition]{Corollary}
 \newtheorem{lemma}[definition]{Lemma}
 \newtheorem{defi}[definition]{Definition}
 
 \def \dem{\noindent{\sc Proof.~}}
\def \findem {\hfill{\hbox {\vrule\vbox{\hrule width 6pt\vskip 6pt\hrule}\vrule}}}

\def\OC{{\mathcal{O}}}
\def\RC{{\mathcal{R}}}
\def\SC{{\mathcal{S}}}

\def\R{{\mathcal{R}}}
\def\g{{\mathfrak{g}}}

\def\a{{\mathfrak{a}}}
\def\b{{\mathfrak{b}}}

\def\m{{\mathfrak{m}}}
\def\mb{{\m_{G^*\times G^*}}}

\def\NM{{\mathbb{N}}}

\def\KM{{\mathbb{K}}}

\def \ho{\widehat{\otimes}}
\def \bo{\bar{\otimes}}

\def \gr {\mathop{\hbox{\rm gr}}\nolimits}
\def \Lift {\mathop{\hbox{\rm Lift}}\nolimits}

\def \Fil {\mathop{\hbox{\rm Fil}}\nolimits}
\def \Prim {\mathop{\hbox{\rm Prim}}\nolimits}

\def \Ker {\mathop{\hbox{\rm Ker}}\nolimits}

\def \lba {\mathop{\underline{\hbox{\rm LBA}}}\nolimits}

\def \Quant {\mathop{\hbox{\rm Quant}}\nolimits}

\def \Braid {\mathop{\hbox{\rm Braid}}\nolimits}
\def \Assoc {\mathop{\hbox{\rm Assoc}}\nolimits}

\def \Aut {\mathop{\hbox{\rm Aut}}\nolimits}
\def \Der {\mathop{\hbox{\rm Der}}\nolimits}

\def \Ad {\mathop{\hbox{\rm Ad}}\nolimits}
\def \Im {\mathop{\hbox{\rm Im}}\nolimits}

\def \cang {\mathop{{\scriptstyle{\rm can}}}\nolimits_{{\g}}}
\def \canget {\mathop{{\hbox{\rm can}}}\nolimits_{{\g^*}}}
\def \cange {\mathop{{\hbox{\rm can}}}\nolimits_{{\g}}}

\def \mod {\mathop{\hbox{\rm mod}}\nolimits}
\def \ad {\mathop{\hbox{\rm ad}}\nolimits}
\def \op {{\scriptstyle{\rm op}}}

\def \dif {{\scriptstyle{\rm d}}}

\def \cohoch {{\scriptstyle{\rm co\hbox{-}Hoch}}}
\def \univ {{\scriptstyle{\rm univ}}}
\def \WX {{\scriptstyle{\rm WX}}}

\def \sym {{\scriptstyle{\rm Sym}}}

\def \univ {{\scriptstyle{\rm univ}}}
\def \Id {\mathop{\hbox{\rm id}}\nolimits}

\def \To {\mathop{\longrightarrow}\limits}
\def \Tto {\mathop{\longleftarrow}\limits}
\def \Hoo {\mathop{\hookleftarrow}\limits}
\def \Simto {\mathop{\sim}\limits}

\def \Uh{U_\hbar(\g)}
\def \Og{\OC_{G^*}}
\def\RW{{\widetilde{\R}_{\WX}}}

\def \dif{{\rm d}}

\begin{document}
\title{On the unicity of braidings of
quasitriangular Lie bialgebras}

\author{Benjamin Enriquez${}^{\dag a}$ Fabio Gavarini${}^\ddag$ 
and Gilles Halbout${}^{\dag b}$\\
{\small{}}\cr
{\small{${}^{\dag}$ \! Institut de Recherche Math\'ematique Avanc\'ee de Strasbourg}}\cr
{\small{UMR 7501 de l'Universit\'e Louis Pasteur et du CNRS}}\cr
{\small{7, rue R. Descartes F-67084 Strasbourg }}\cr
{\small{${}^{\dag a}$ \! e-mail:\,\texttt{enriquez@math.u-strasbg.fr}}}\cr
{\small{${}^{\dag b}$ \! e-mail:\,\texttt{halbout@math.u-strasbg.fr}}}\cr
{\small{}}\cr
{\small{${}^{\ddag}$ \! Dipartimento di Matematica}}\cr
{\small{Universit\`a degli Studi di Roma ``Tor Vergata''}}\cr
{\small{Via della Ricerca Scientifica, 1   ---   I-00133 Roma, Italy}}\cr
{\small{\! e-mail:\,\texttt{gavarini@mat.uniroma2.it}}}}
\markboth{Benjamin Enriquez Fabio Gavarini and Gilles Halbout}
{On the unicity of braidings of
quasitriangular Lie bialgebras}

\maketitle

\abstract{Any quantization of a quasitriangular Lie bialgebra $\g$ gives rise
to a braiding of the dual Poisson-Lie formal group $G^*$.
We show that this braiding always coincides with the Weinstein-Xu braiding.
We also define lifts of the classical $r$-matrix a certain formal 
functions on $G^* \times G^*$, prove their existence
and uniqueness using co-Hochschild cohomology arguments and show 
show that the lift can be expressed in terms of $r$ by universal formulas.}

\vskip20pt

\centerline {\bf \S\; 0 \ Outline of results}

\vskip20pt

\bigskip

\noindent {\bf - a - Quasitriangular Lie algebras}

\bigskip

We fix a base field $\KM$ of characteristic zero. Let $(\g,r)$ be a 
finite dimensional
quasitriangular Lie bialgebra. Recall that this means that
\begin{itemize}
\item $(\g,[-,-],\delta)$ is a Lie bialgebra;
\item $r \in \g \otimes \g$ is a solution of the classical Yang-Baxter equation
(CYBE), i.e., 
$$[r^{1,2},r^{1,3}]+[r^{1,2},r^{2,3}]+[r^{1,3},r^{2,3}]=0;$$
\item we have $\delta(x)=[r,x\otimes 1+1\otimes x]$ for any $x \in \g$,
so in particular,
$r+r^{2,1}$ is $\g$-invariant.
\end{itemize}

\bigskip

\noindent {\bf - b - $\Quant(\g)$}

\bigskip

A {\it quantization} of $(\g,r)$ is a quantized universal enveloping (QUE) algebra
$(\Uh,m,\Delta)$ quantizing $(\g,[-,-],\delta)$, together with an element $R \in
\Uh^{\ho 2}$, such that if $x \mapsto (x \mod \hbar)$ is the canonical
projection
$\Uh^{\ho 2} \to U(\g)^{\ho 2}$ then
\begin{itemize}
\item $\Delta^\op=R\Delta R^{-1},$
\item $(\Delta \otimes \Id)(R)=R^{1,3}R^{2,3}$, 
$(\Id \otimes \Delta)(R)=R^{1,3}R^{1,2}$,
\item $(\epsilon \otimes \Id)(R)=(\Id \otimes \epsilon)(R)=1$
where $\epsilon$~: $\Uh \to \KM[[\hbar]]$ is the counit of $\Uh$,
\item $(R \mod \hbar)=1$, $\left({{R-1}\over{\hbar}} \mod \hbar \right)=r \in \g
\otimes \g \subset U(\g) \otimes U(\g)$.
\end{itemize}
\noindent We denote by $\Quant(\g)$ the set of all quantizations
of $(\g,r)$. According to \cite{EK}, we have a map
$\Assoc(\KM) \to \Quant(\g)$
(where $\Assoc(\KM)$ is the set of all Lie associators defined over $\KM$), so $\Quant(\g)$ is nonempty.

\bigskip

\noindent {\bf - c - $\Braid(\g)$}

\bigskip

Let $G^*$ be the formal group corresponding to the dual Lie bialgebra $\g^*$,
and let $\Og$ be its function ring; so $\Og=(U(\g^*))^*$; this is a formal
series Hopf algebra, equiped
with coproduct 
$\Delta_\OC$~: $\Og \to \Og \bo \Og$ ($\bo$ is the tensor product
of the formal series algebras, $\Og \bo \Og=\OC_{G^* \times G^*}$ 
is the function ring of $G^* \times G^*$)
and
counit $\epsilon_\OC$~: $\Og \to \KM$.

\begin{defi}
\label{def:0.1}
A braiding of $G^*$ is a Poisson
algebra automorphism $\RC$ of  
$\Og \bo \Og$ satisfying the conditions:
\begin{description}
\item[($\alpha$)] $(\epsilon_\OC \otimes \Id)\circ \RC=\epsilon_\OC \otimes \Id,~
(\Id \otimes \epsilon_\OC)\circ \RC=\Id \otimes \epsilon_\OC,$
\item[($\beta$)] $\Delta_\OC^{\op}=\RC \circ \Delta_\OC,$
\item[($\gamma$)] $\RC^{1,3} \circ \RC^{2,3}\circ (\Delta_\OC \otimes \Id)=
(\Delta_\OC \otimes \Id) \circ \RC,$
\newline
$\RC^{1,3} \circ \RC^{1,2}\circ (\Id \otimes \Delta_\OC)=
(\Id \otimes \Delta_\OC) \circ \RC,$
\item[($\delta$)] if $\mb$ is the maximal ideal of $\Og \bo \Og$, then
\begin{itemize}
\item the automorphism $\m_{G^* \times G^*}/\m_{G^* \times G^*}^2
\to \m_{G^* \times G^*}/\m_{G^* \times G^*}^2$ induced by $\R$ is the identity, 
\item 
therefore $\R -\Id$ induces a linear map
$[\R-\Id]~:~\m_{G^* \times G^*}/\m_{G^* \times G^*}^2$ $ \to 
\m_{G^* \times G^*}^2/\m_{G^* \times G^*}^3,$
and if we use the natural identifications 
\begin{align*}
\m_{G^* \times G^*}/\m_{G^* \times G^*}^2 &\To^\sim \g \oplus \g \cr
\m_{G^* \times G^*}^2 / \m_{G^* \times G^*}^3 &\To^\sim
\left(S^2(\g)\right) \oplus \left( \g \otimes \g \right)
\oplus \left( S^2(\g)\right),
\end{align*}
then
$[\R -\Id]$ coincides with the map
$$(x,y)\mapsto (0,[r,x\otimes1+1\otimes y],0).$$ 
\end{itemize}
\end{description}
\end{defi}
\noindent We denote by $\Braid(\g)$ the set of all braidings of $G^*$.

\bigskip

\noindent {\bf - d - The Weinstein-Xu braiding}

\bigskip

Define $\RW$~: $G^* \times G^* \to G^* \times G^*$ by 
\begin{equation}
\label{Req}
\RW(u,v)=(\lambda_{R_-(v)}(u),\rho_{R_+(u)}(v)),
\end{equation}
where $R_\pm$~: $G^* \to G$ are the formal group morphisms exponentiating
the Lie algebra morphisms
$r_\pm$~: $\g^* \to \g$, where 
$r_+(\xi)=\langle r,\xi\otimes \Id \rangle$ and
$r_-(\xi)=-\langle r,\Id \otimes \xi\rangle $, and $\lambda$, $\rho$ are the left and right
dressing actions of $G$ on $G^*$ (regular action on $G^*=D/G$
and on $G^*=D\backslash G$, where $D$ is the double group of $G$).

\smallskip

Let $\R_{\WX} \in \Aut(\Og \bo \Og)$ be the algebra automorphism induced be
$\RW$. Then
$$\R_{\WX} \in \Braid(\g)~\hbox{(see \cite{WX} and \cite{GH2})}.$$

\bigskip

\noindent {\bf - e - The Gavarini-Halbout map}

\bigskip

If $(\Uh(\g),m,\Delta,R)$ is a quantization of $(\g,r)$, define $\OC_\hbar$ as a
quantized function algebra associated to $\Uh$.
So
$$\OC_\hbar=\{f \in \Uh|~\forall n \geq 0,~\delta^{(n)}(f)
\in \hbar^n \Uh^{\ho n}\};$$
where $\delta^{(n)}$~: $\Uh \to \Uh^{\ho n}$ is defined by
$\delta^{(n)}=(\Id -\eta \circ \epsilon)^{\otimes n} \circ \Delta^{(n)}$.
Then $\OC_\hbar$ is a topological Hopf subalgebra of $\Uh$, and it is a
quantization of the Hopf-Poisson algebra $\Og$ (see [Dr,Ga]).
In particular, $\OC_\hbar/\hbar\OC_\hbar\simeq \Og$.
\begin{theorem}
\label{thm:GH}
(see \cite{GH} and also \cite{EH}) The inner automorphism
$\Ad(R)$~: $x \mapsto RxR^{-1}$ of $\Uh^{\ho 2}$ restricts to an
automorphism $\RC_\hbar$ of $\OC_\hbar^{\bo 2}$.
The reduction $\R$ of $\R_\hbar$ modulo $\hbar$ is
an outer automorphism of $\Og \bo \Og$. Then $\RC \in \Braid(\g)$. 
\end{theorem}
\noindent 
The main part of this result was proved in \cite{GH} (see also \cite{EH}).
The remaining part is a consequence of Proposition \ref{prop:0.7}.
Therefore we have a map:
$$\hbox{GH~: }\Quant(\g) \to \Braid(\g).$$

\bigskip

\noindent {\bf - f - Unicity of braidings}

\bigskip

\begin{theorem}
\label{thm:braidings}
\label{thm:3}
$\Braid(\g)$ contains only one element, so
$$\Braid(\g)=\{\R_{\WX}\}.$$
\end{theorem}

\medskip

\noindent In particular, the braiding $\R$ constructed in Theorem \ref{thm:GH}
coincides with $\R_{\WX}$.

\bigskip

\noindent {\bf - g - Formal Poisson manifolds}

\bigskip

Let $A$ be an arbitrary Poisson formal series algebra; let us denote by $\m_A$
the
maximal ideal of $A$, and let us assume that $\{A,A\}\subset \m_A$.
Then we have
$\{\m_A^k,\m_A^l\} \subset \m_A^{k+l-1}$, for any $k,l \geq 0$.
For $f,g \in \m_A^2$, the
Campbell-Baker-Hausdorff (CBH) series
$$f \star g = f+g +{{1}\over{2}} \{f,g\} + \cdots + B_k(f,g)+ \cdots $$
converges in $A$.

\smallskip

There is a unique Lie algebra morphism 
\begin{align*}
V~:~A &\to \Der(A)\\
f &\mapsto (V_f~:~g \mapsto \{f,g\}).
\end{align*}
Define $\Der^+(A)$ as the Lie subalgebra of $\Der(A)$ of all derivations taking
each 
$\m_A^k$ to $\m_A^{k+1}$. Then $V$ restricts to a Lie algebra morphism
$\m_A^2 \to \Der^+(A)$. Moreover, for any derivation
$D \in \Der^+(A)$, the series $\exp (D)$ is a well defined automorphism of $A$;
this defines an exponential map
\begin{align*}
\exp\! ~:~\Der^+(A) &\to \Aut(A)\\
D &\mapsto \exp\! (D).
\end{align*}
The series $\exp(D)$ is a well-defined automorphism of $A$.
Let us denote by $\Aut^+(A)$ the subgroup of $\Aut(A)$ of all Poisson
automorphisms $\theta$ such that the
map $[\theta]$~: $\m_A/\m_A^2 \to
\m_A/\m_A^2$ induced by $\theta$ is the identity.
Then $\exp(D)$ belongs to $\Aut^+(A)$, and the map
$\exp$~: $\Der^+(A) \to \Aut^+(A)$ is a bijection.

\bigskip

\noindent {\bf - h - Lifts of the classical $r$-matrix}

\bigskip

Using the previous section for the formal Poisson manifold $\Og$, 
we can define lifts:
\begin{defi}
\label{def:lifts}
\label{def:0.4}
A lift of $r$ is an element $\rho \in \Og \bo \Og$, such that:
\begin{description}
\item[($\alpha$)] $(\epsilon \otimes \Id)(\rho)=
(\Id \otimes \epsilon)(\rho)=0$,
\item[($\beta$)] $\Delta^{\op}=\Ad(\exp (V_\rho))\circ \Delta$
(equality of automorphisms of 
$\OC_{G^* \times G^*}$),
\item[($\gamma$)] $(\Delta \otimes \Id)(\rho)=\rho^{1,3}\star
\rho^{2,3},\hskip0.3cm
(\Id \otimes \Delta)(\rho)=\rho^{1,3}\star \rho^{1,2}$, 
where $\rho^{i,j}$ is the image of $\rho$ by
the
map $\left(\OC_{G^*}\right)^{\bo 2} \to \left(\OC_{G^*}\right)^{\bo 3}$
associated with $(i,j)$,
\item[($\delta$)] the class $[\rho]$ of $\rho$ in
$\left(\m_{G^*}/\m_{G^*}^2\right)^{\otimes 2}=\g \otimes \g$ satisfies
$$[\rho]=r.$$
\end{description}
\end{defi}
\noindent 
Condition ($\beta$) may be rewritten as follows:
$$\forall f \in \Og,~\Delta^\op(f)= \rho \star \Delta(f) \star(-\rho).$$
It will follow from the proof of Theorem \ref{thm:0.8}
that this condition may be dropped
from the definition of
$\Lift(\g)$ (see Lemma \ref{lift:lift'}).
We denote by $\Lift(\g)$ the set of all lifts of $r$.

\bigskip

\noindent {\bf - i - Sequence of maps ${\Quant(\g)\! \to \! \Lift(\g)\!  \to
\! \Braid(\g)}$}

\bigskip

Let us recall an $\hbar$-adic valuation result for $R$-matrices:

\begin{theorem}{(\cite{EH})}
If $(\Uh,m,\Delta,R)$ is a quantization of $(\g,r)$, and if we set
$\rho_\hbar=\hbar\log(R)$, then $\rho_\hbar \in \OC_\hbar^{\bo 2}$.
If $\m_\hbar $ is the kernel of the counit map
$\OC_\hbar \to \KM[[\hbar]]$, we
even have $\rho_\hbar \in \m_\hbar^{\bo 2}$.
\end{theorem}

\begin{corollary}
\label{cor:O.6}
The reduction $\rho$ of $\rho_\hbar$ modulo $\hbar$ belongs to $\Lift(\g)$. So
the assignment
$(\Uh,m,\Delta,R) \mapsto (\rho_\hbar \mod \hbar)$ defines a map
$\Quant(\g) \to \Lift(\g)$.
\end{corollary}

\begin{proposition}
\label{prop:0.7}
There is a unique map
$\Lift(\g) \to \Braid(\g)$, taking $\rho$ to
$\exp\! (V_\rho)$. Then the composed map
$\Quant(\g) \to \Lift(\g) \to \Braid(\g)$ coincides with
{\rm GH}~: $\Quant(\g) \to \Braid(\g)$.
\end{proposition}

\bigskip

\noindent {\bf - j - Unicity of lifts}

\bigskip

\begin{theorem}
\label{thm:lift}
\label{thm:0.8}
$\Lift(\g)$ consists of only one element.
\end{theorem}
\noindent The unicity part of this theorem uses an elementary argument.
The existence part uses the nonemptiness of
$\Quant(\g)$, so it relies on the theory of associators and transcendental
arguments. In the last part of the paper, we outline
an algebraic proof of the existence part of Theorem \ref{thm:lift}, relying on
co-Hochschild cohomology arguments.

\vskip1cm 

\noindent {\bf - k - Universal versions}

\bigskip

If $\a$ is a finite dimensional Lie bialgebra and $\g$ is the double
of $\a$ (so $\g=\a \oplus \b$, $\b=\a^*$), then we have the algebra
isomorphisms
$$\Og\simeq \widehat{S}^\cdot(\g) \simeq \widehat{S}^\cdot(\a)\bo
\widehat{S}^\cdot(\b)$$
where $\widehat{S}^\cdot$ is the graded completion of the symmetric algebra.
The last isomorphism is dual to the composed map
$${S}^\cdot(\b)\otimes S^\cdot(\a) \To^{\sym \otimes \sym} U(\a) \otimes
U(\b) \To^m U(\g)$$
where Sym is the symetrization map.

\medskip

Therefore $\Og^{\bo n} \simeq \widehat{S}^\cdot(\a)^{\bo n} \bo
\widehat{S}^\cdot(\b)^{\bo n}.$
Now if $F$ and $G$ are any Schur functors, one can define a universal
version
of the space $F(\a) \otimes G(\b)$, namely
$(F(\a)\! \otimes G(\b))_{\univ}\!=\!\lba(G,F)$, where
$\lba$ is the prop of Lie bialgebras (see, e.g.,  \cite{EE}).
We then define Hopf algebras $\left(\Og^{\bo n}\right)_{\univ}=
\left(\widehat{S}^\cdot(\a)^{\bo n} \bo
\widehat{S}^\cdot(\b)^{\bo n}\right)_\univ$,
together with insertion-coproduct morphisms relating them.

\begin{defi}
A universal lift is an element $\rho_\univ \in
\left(\Og^{\bo 2}\right)_{\univ}$, satisfying the universal
versions of the conditions of Definition \ref{def:lifts}.
\end{defi}
\noindent We denote by $\Lift_\univ$ the set of all universal lifts.

\medskip

When $\g$ is any finite-dimensional quasitriangular Lie bialgebra, we have
algebra morphisms
$\left(\Og^{\bo n}\right)_{\univ} \to \Og^{\bo n}$.
It follows that for any $\g$, we have a map
$\Lift_\univ \to \Lift(\g)$.
\begin{theorem}
\label{thm:1.10}
$\Lift_\univ$ consists of only one element $\rho_\univ$.
\end{theorem}
\noindent So the unique lift $\rho_\g$ of a quasitriangular Lie bialgebra
$\g$ is obtained from the element
$$r \in \g \otimes \g \subset \widehat{S}^\cdot(\g) \bo \widehat{S}^\cdot(\g)
=\Og^{\bo 2}$$
by universal formulas.
In \cite{Re}, Reshetikin computed $\rho_\g$ when
$\g$ is a semi-simple Lie algebra. His formulas
involve the dilogarithm function. We do not know an explicit formula for
$\rho_\univ$.
It might be simpler to express the pairing
$\langle -,-\rangle $~: $U(\g^*)^{\otimes 2}\to \KM$, defined by
$\langle x,y\rangle =\langle \rho_\g,x\otimes y\rangle ;$ this way one avoids 
the unnatural use of
symmetrization maps.

\bigskip

\noindent {\bf - l - Plan of the paper}

\bigskip

\noindent In Section 1, we construct a map $\Quant(\g)\to \Lift(\g)$
(Corollary \ref{cor:O.6}) and prove the unicity of lifts
(Theorem \ref{thm:0.8}).

\smallskip

\noindent In Section 2, we construct the map $\Quant(\g) \to \Braid(\g)$
(Proposition \ref{prop:0.7}), and then prove the unicity of braidings (Theorem
\ref{thm:3}).
The proof of this theorem uses only a part of the arguments of Section 1
(essentially only the existence of a sequence of maps
$\Quant(\g)\to \Lift(\g) \to \Braid(\g)$).

\smallskip

\noindent In Section 3, we outline a proof of Theorem
\ref{thm:0.8} not depending on the theory of associators. 

\smallskip

\noindent In Section 4, we sketch a proof of Theorem \ref{thm:1.10}

\smallskip

\noindent In Section 5 (appendix), we construct a commutative diagram related to the
duality theory of quantized universal enveloping algebras, which we use in the
Sections 1 and 2.

\vskip20pt

\centerline {\bf \S\; 1 \ Lifts of classical $r$-matrices}

\stepcounter{section}\label{section:1}

\vskip20pt

\begin{proposition}
There exists a map $\Quant(\g) \to \Lift(\g)$.
\label{prop:existence}
\end{proposition}

\dem
Let $(U_\hbar(\g),m, \Delta_\hbar)$ be an element of $\Quant(\g)$.
Let $\OC_\hbar\subset U_\hbar(\g)$ be the quantized formal series Hopf
(QFSH)
subalgebra sitting in $U_\hbar(\g)$. Let $\m_\hbar$ be the augmentation
ideal of $\OC_\hbar$;
then $\m_\hbar\subset \hbar U_\hbar(\g)$.
In \cite{EH}, we showed that there exists a unique $\rho_\hbar \in
\m_\hbar^{\bo 2}$ such that $R=\exp\! \left({{\rho_\hbar}\over{\hbar}}\right)$
(this exponential is well-defined because
$\left({{\rho_\hbar}\over{\hbar}}\right) \in \hbar U_\hbar(\g)^{\ho 2}$).
Then the quasitriangular identities of $R$ can be translated as follows: for
$a,b \in \OC_\hbar^{\bo 3}$, we  set
$\{a,b\}_\hbar = {{1}\over{\hbar}}[a,b]$. Let
$\m_\hbar^{(3)}$ be the augmentation ideal of
$\OC_\hbar^{\bo 3}$. Then
$$\{\m_\hbar^{(3)},\m_\hbar^{(3)}\}_\hbar \subset \m_\hbar^{(3)},$$
therefore
$\left\{\left(\m_\hbar^{(3)}\right)^k,\left(\m_\hbar^{(3)}\right)^l
\right\}_\hbar \subset \left(\m_\hbar^{(3)}\right)^{k+l-1}.$
Now if $a,b \in \left(\m_\hbar^{(3)}\right)^2$, the series
$$a \star_\hbar b = a +b + {{1}\over{2}} \{a,b\}_\hbar + \cdots $$
(CBH series, where the Lie bracket is $\{-,-\}_\hbar$) is convergent
in $\OC_\hbar^{\bo 3}$. Then we have~:
\begin{equation}
(\Delta_\hbar \otimes \Id)(\rho_\hbar)=\rho_\hbar^{1,3}\star_\hbar
\rho_\hbar^{2,3},\hskip0.5cm
(\Id \otimes \Delta_\hbar)(\rho_\hbar)=\rho_\hbar^{1,3}\star_\hbar 
\rho_\hbar^{1,2}.
\label{interm}
\end{equation}
$\Delta_\hbar$ restricts to a map $\OC_\hbar \to \OC_\hbar^{\bo 2}$,
the reduction of which modulo $\hbar$ is the coproduct map $\Delta$ of
$\OC_{G^*}$. 
Define $\rho$ as the reduction modulo $\hbar$ of $\rho_\hbar$, so 
$\rho \in \m_{G^*}^{\bo 2}$.
Taking the reduction of (\ref{interm}) modulo $\hbar$, we get
($\gamma$) of Definition \ref{def:0.4}.

\smallskip

On the other hand, we have
$\Delta_\hbar^{\op}=\Ad\left(\exp\! \left({{\rho_\hbar}\over{\hbar}}\right)\right)\circ
\Delta$. Set $\ad_\hbar(a)(b)=\{a,b\}_\hbar$. The automorphisms $\Ad
\left(\exp\! \left({{\rho_\hbar}\over{\hbar}}\right)\right)$
and $\exp\! (\ad_\hbar(\rho_\hbar))$ coincide.
So we get the identity~:
\begin{equation}
\Delta_\hbar^{\op}=\exp\! (\ad_\hbar(\rho_\hbar))\circ \Delta_\hbar
\label{interm:2}
\end{equation}
(equality of two morphisms $\OC_\hbar\to \OC_\hbar^{\bo
  2}$).
Taking the reduction of (\ref{interm:2}) modulo $\hbar$, we get
($\beta$) of Definition \ref{def:0.4}.

\smallskip

To show that $\rho$ satisfies ($\delta$) of Definition \ref{def:0.4}, we use the following
result (which will be proved in Section 4)~:
\begin{lemma}
Let $\sigma$  be an arbitrary element of $\m_\hbar \bo \m_\hbar$
and $[\sigma]$ be its class in
$\left(\m_\hbar/\hbar \m_\hbar + \m_\hbar^2\right)^{\bo 2}$.
Since $\left(\m_\hbar/\hbar \m_\hbar + \m_\hbar^2\right)$ identifies
with $\g$,
$[\sigma]\in \g^{\otimes 2}$.
Since $\m_\hbar \subset \hbar U_\hbar(\g)$,
$\sigma$ is an element of  $\hbar^2 U_\hbar(\g)^{\ho 2}$.
Then $\left({{\sigma}\over{\hbar^2}} \mod \hbar\right)$ is an element of
$U(\g)^{\otimes 2}$.
We have the following identity in $U(\g)^{\otimes 2}$~:
$$\left({{\sigma}\over{\hbar^2}}\mod \hbar \right)=[\sigma].$$
\label{gav}
\end{lemma}
\noindent Then $\m_\hbar \subset \hbar U_\hbar(\g)$, so
${{\rho_\hbar}\over{ \hbar}} \in \hbar U_\hbar(\g)^{\ho 2}$, so
$${{R-1}\over{\hbar}}={{\rho}\over{\hbar^2}}+{{1}\over{2\hbar}}\left(
{{\rho}\over{\hbar}}\right)^2+{{1}\over{6\hbar}}\left({{\rho}\over{\hbar}}\right)^3+\cdots
$$
The terms 
${{1}\over{2\hbar}}\left(
{{\rho}\over{\hbar}}\right)^2$,
${{1}\over{6\hbar}}\left({{\rho}\over{\hbar}}\right)^3$,$\dots$, all
belong
to $\hbar U_\hbar(\g)^{\ho 2}$ so
\begin{equation}
\left({{R-1}\over{\hbar}} \mod \hbar\right)=\left({{\rho}\over{\hbar^2}}
    \mod \hbar\right).
\label{id:rho}
\end{equation} 
Now
\begin{align*}
[\rho]&=\left({{\rho}\over{\hbar^2}} \mod \hbar\right)&\hbox{(by Lemma
  \ref{gav})}\\
&=\left({{R-1}\over{\hbar}} \mod \hbar\right)&\hbox{(by identity
  (\ref{id:rho}))}\\
&=r &\hbox{(by hypothesis on }R).
\end{align*}
Therefore $\rho$ satisfies property ($\delta$) of Definition \ref{def:0.4}.

\smallskip

\noindent Now we have proved that the reduction $\rho$ of 
$\rho_\hbar$ modulo $\hbar$ satisfies all
the conditions of
Definition \ref{def:lifts}.
\findem

\medskip

\begin{proposition}
$\Lift(\g)$ contains at most one element.
\label{prop:unicity}
\end{proposition}

\dem
Let us denote by $\m_{G^* \times G^*}$ the maximal ideal of
$\OC_{G^* \times G^*}$, so
$\m_{G^* \times G^*}$ $=\m_{G^*}\bo \OC_{G^*}+\OC_{G^*}\bo \m_{G*}$.
Then we have for any $N \geq 0$, 
$$\m_{G^*}^{\bo 2}\cap \m_{G^* \times G^*}^N=\sum_{
{\begin{array}{l}
{\scriptstyle{a,b \geq 1}}\cr
{\scriptstyle{a+b=N}}
\end{array}}
}\m_{G^*}^a \bo \m_{G^*}^b.$$
Let $\rho$ and $\rho'$ be two lifts or $r$. The classes of $\rho$ and $\rho'$ are the same in
$\m_{G^*}^{\bo 2}/(\m_{G^*}^{\bo 2} \cap \m_{G^* \times G^*}^2)$ and
  equal to $r$, by assumption.

\smallskip

Let $N$ be an integer $\geq 2$; assume that we have proved that $\rho$
and
$\rho'$ are equal modulo 
$\m_{G^*}^{\bo 2} \cap \m_{G^* \times G^*}^N$.
Let us show that they are equal modulo
$\m_{G^*}^{\bo 2} \cap \m_{G^* \times G^*}^{N+1}$.
Write $\rho'=\rho+\sigma$; then $\sigma \in \m_{G^*}^{\bo 2}\cap
\m_{G^* \times G^*}^N$.
We get

\begin{align}
\label{id:sigma}
\begin{split}
\left(\Delta \otimes \Id\right)(\sigma)=&
(\rho + \sigma)^{1,3}\star (\rho + \sigma)^{2,3}-\rho^{1,3} \star
\rho^{2,3}\\
=& \sigma^{1,3} + \sigma^{2,3}\\
&+\sum_{k>1}\left(B_k\left( \rho^{1,3} + \sigma^{1,3},
\rho^{2,3} + \sigma^{2,3}\right) - B_k\left(\rho^{1,3},\rho^{2,3}\right)\right),
\end{split}
\end{align}
where $B_k$ is the total degree $k$ Lie polynomial of the CBH series.

\begin{lemma}
If $k>1$,
$B_k\left(\rho^{1,3}+\sigma^{1,3},\rho^{2,3}+\sigma^{2,3}\right)
-B_k\left(\rho^{1,3},\rho^{2,3}\right)$ is an element of $ \m_{G^*\times G^*}^{N+1}$.
\label{lemma:approx}
\end{lemma}

\dem
This difference may be expressed as a sum of terms of the form
$$P_k\left(\sigma^{i_1,3},\dots,\sigma^{i_l,3},\rho^{i_{l+1},3},\dots,\rho^{i_k,3}\right),$$
where $P_k$ is a Lie polynomial, homogeneous of degree $1$ in
each variable $i_1,\dots,i_k \in \{1,2\}$, and
$l \geq 1$. This expression belongs to
$\m_{G^*\times G^*}^{l(N-2)+k+1}$.
So it belongs to
$\m_{G^* \times G^*}^{N+k-1} \subset \m_{G^* \times G^*}^{N+1}$.
\findem

\smallskip

Now $\OC_{G^*}$ is equipped with a decreasing Hopf filtration
$\OC_{G^*}\supset \m_{G^*} \supset \m_{G^*}^2 \supset \cdots$~:
we have
$$\Delta\left(\m_{G^*}^k\right)\subset \sum_{\alpha,\beta|\alpha+\beta =k}
\m_{G^*}^\alpha \bo \m_{G^*}^\beta.$$
Its associated graded is therefore also a Hopf algebra; it is
isomorphic to the formal completion $\widehat{S}^\cdot(\g)$ of the
commutative
and cocommutative symmetric algebra $S^\cdot(\g)$, the coproduct of which is
defined by the condition
that the elements of degree $1$ are
primitive.
The tensor square $\OC_{G^*}^{\bo 2}$ is also
filtered: the $i$-th term of the decreasing filtration is
$$\Fil^i\left(\OC_{G^*}^{\bo 2}\right)=\sum_{\alpha,\beta|\alpha +\beta =i}
\m_{G^*}^\alpha \bo \m_{G^*}^\beta;$$
and we have
$$\gr\left(\OC_{G^*}^{\bo 2}\right)=\widehat{S}^\cdot(\g) \bo \widehat{S}^\cdot(\g).$$
Moreover, let $[\sigma]$ be the class of $\sigma$ in 
$\gr^N\left(\OC_{G^*}^{\bo 2}\right)$;
according to identity (\ref{id:sigma}) and 
Lemma \ref{lemma:approx}, we have
$$\left(\Delta \otimes \Id\right)([\sigma])=
[\sigma]^{1,3}+[\sigma]^{2,3},\hskip0.5cm
\left(\Id \otimes \Delta \right) ([\sigma])=[\sigma]^{1,3} + [\sigma]^{1,2}.$$
The first identity implies that $[\sigma]\in \g \otimes S^{N-1}(\g)$, the second
identity implies that $[\sigma]\in S^{N-1}(\g) \otimes \g$; since
$\left(\g \otimes S^{N-1}(\g)\right) \cap \left(S^{N-1}(\g)\otimes \g 
\right)=\{0\}$, we get $[\sigma]=0$, therefore $\sigma \in
\m_{G^* \times G^*}^{N+1}.$
So $\sigma $ belongs to the intersection of
all $\m_{G^* \times G^*}^N$, 
$N \geq 0$, so $\sigma=0$. This proves that $\rho = \rho'$.
\findem

\smallskip

\begin{corollary}
\label{cor:ex:unicity}
If $(\g,r)$ is a quasitriangular Lie bialgebra, there exists a unique element
$\rho \in \Lift(\g)$.
\end{corollary}

\dem
The unicity follows from Proposition \ref{prop:unicity},
and the existence follows from Proposition \ref{prop:existence},
and from the fact that $\Quant(\g)$ is nonempty:
in \cite{EK}, Etingof and Kazhdan constructed a map
$\Assoc(\KM)\to \Quant(\g)$, where $\Assoc(\KM)$ is the set of
associators over the ground field $\KM$;
this set is introduced by Drinfeld in \cite{Dr},
where it is also shown that $\Assoc(\KM)$ is nonempty.
\findem

\begin{remark}
Corollary \ref{cor:ex:unicity} relies on the existence of associators, so it
actually relies on transcendental arguments. Another
proof of this
Corollarary will be given in Section 3;
this proof is algebraic and is based on the further use of
co-Hochschild cohomology groups.
\end{remark}

\vskip20pt

\centerline {\bf \S\; 2 \ Quasitriangular braidings}

\stepcounter{section}\label{section:2}

\vskip20pt

In this section, we construct the map
$\Quant(\g) \to \Braid(\g)$ (Subsection 2.a).
We then prove that 
$\RC_{\WX} \in \Braid(\g)$ (Subsection 2.b).
In Subsection 2.c, we prove that $\Braid(\g)$ contains at most one
element.
So the image of any element of $\Quant(\g)$ in
$\Braid(\g)$ coincides with $\RC_{\WX}$ (Theorem \ref{thm:3}).

\bigskip

\noindent {\bf - a - The map $\Quant(\g)\to \Braid(\g)$ (proof of Proposition
\ref{prop:0.7})}

\bigskip

Let us prove the map $\rho \mapsto \exp\! \left(V_\rho\right)$ actually maps
$\Lift (\g) \to \Braid(\g)$.
If $\rho \in \Lift(\g)$, the
fact that $\rho$ satisfies axioms ($\alpha$), ($\beta$) and ($\gamma$)
of Definition \ref{def:0.4} respectively implies that 
$\exp\! \left(V_\rho\right)$ satisies axioms ($\alpha$), ($\beta$) and ($\gamma$)
of Definition \ref{def:0.1}. Let us now prove that the fact that $\rho$ 
satisfies axiom ($\delta$) of Definition \ref{def:0.4} implies
that $\exp\! \left(V_\rho\right)$ satisfies axiom ($\delta$)
of Definition \ref{def:0.1}.

\medskip

By definition, $\rho$ is an element of  $\m_{G^*}\bo \m_{G^*}$.
We have 
$\m_{G^*}\bo \m_{G^*} \subset 
\OC_{G^*}\bo \OC_{G^*}=\OC_{G^*\times G^*}$; actually, we have
$\m_{G^*}\bo \m_{G^*} \subset \m_{G^*\times G^*}^2$,
so $\rho \in \m_{G^*\times G^*}^2$.
Since we have 
$\{\m_{G^* \times G^*}^2,\m_{G^* \times G^*}\}
\subset \m_{G^* \times G^*}^2$, the map 
$$V_\rho~:~\m_{G^* \times G^*}/\m_{G^* \times G^*}^2\to
\m_{G^* \times G^*}/\m_{G^* \times G^*}^2$$ 
induces the zero map.
Therefore, so do all
the
$\left(V_\rho\right)^k$, $k\geq 1$.
So $\exp\! \left(V_\rho\right)$ induces the identity map
of $\m_{G^* \times G^*}/\m_{G^* \times G^*}^2$.
Let us now compute the map
$$[\exp\! \left(V_\rho\right)-\Id]~:~
\m_{G^* \times G^*}/\m_{G^* \times G^*}^2
\to \m_{G^* \times G^*}^2 /\m_{G^* \times G^*}^3$$
using the identifications 
$ \m_{G^* \times G^*}^2 /\m_{G^* \times G^*}^3=
S^2(\g)\oplus (\g \otimes \g) \oplus S^2(\g)$ and
$\m_{G^* \times G^*}/\m_{G^* \times
G^*}^2=\g \oplus \g$.
We have 
$$(\epsilon \otimes \Id) \circ \left(\exp\! \left(V_\rho\right)-\Id\right)=
(\Id \otimes \epsilon)\circ \left(\exp\! \left(V_\rho\right)-\Id\right)=0$$
(identity of maps $\OC_{G^* \times G^*} \to \OC_{G^*}$), because
$\{\m_{G^*},\OC_{G^*}\}\subset \m_{G^*}$. The class of
$\Ker(\epsilon \otimes \Id)\cap \Ker (\Id \otimes \epsilon) \cap 
\m_{G^* \times G^*}^2$ in $\m_{G^* \times G^*}^2/\m_{G^* \times G^*}^3$
is the subspace
$(\g \otimes \g) \subset S^2(\g) \oplus (\g \otimes \g) \oplus S^2(\g)$.

\smallskip

On the other hand, if $f \in \m_{G^*\times G^*}$, and
$k \geq 2$, then $\left(V_\rho\right)^k (f) \in \m_{G^*\times
G^*}^3$. So the class of $\left(\exp\!
\left(V_\rho\right)-\Id\right)(f)$ in 
$\m_{G^*\times G^*}^2/\m_{G^*\times
G^*}^3$ coincides with that of $V_\rho(f)$.
So we now compute the map
$$[V_\rho]~:~\g \oplus \g \to (\g \otimes \g).$$
Let $ x_1,x_2 \in \g$ and let $f_1,f_2 \in \m_{G^*}$ be such that their classes
in $\m_{G^*}/\m_{G^*}^2=\g$ are $x_1,x_2$. Let us
set 
$f=f_1 \otimes 1 + 1 \otimes f_2$, and let us compute
$V_\rho(f)$.
Set $\rho=\sum_\alpha \rho_\alpha' \otimes
\rho_\alpha''$,
with $\rho_\alpha',\rho_\alpha'' \in \m_{G^*}$.
Then
$$V_\rho(f)=\sum_\alpha\{\rho_\alpha',f_1\}\otimes \rho_\alpha''
+\rho_\alpha'\otimes\{\rho_\alpha'',f_2\}.$$
Now we have a commutative diagram
$$\begin{matrix}
{\m_{G^*}\otimes \m_{G^*}}
&\xrightarrow{\hbox{Poisson bracket}}&\m_{G^*}\cr
\downarrow&&\downarrow\cr
{\g \otimes \g}
&\xrightarrow{\hbox{\hskip0.4cm Lie bracket\hskip0.4cm}}
&{\g}
\end{matrix}$$
where the vertical arrows correspond to the projection $\m_{G^*}\to
\m_{G^*}/\m_{G^*}^2=\g$.
So the class of $V_\rho(f)$ in $\g \otimes \g$ is
$[r,x_1 \otimes 1 + 1 \otimes x_2]$.
Therefore $[\exp\! \left(V_\rho\right)-\Id]~:~\g \oplus \g \to (\g \otimes \g)$
is the map
$$(x_1,x_2)\mapsto\left(0,[r,x_1 \otimes1 + 1 \otimes x_2],0\right),$$
which proves that
$\exp\! \left(V_\rho\right)$ satisfies condition ($\delta$)
of Definition \ref{def:0.1} and so belongs to $\Braid(\g)$.
\findem

\bigskip

\noindent {\bf - b - Proof of $\RC_{\WX}\in \Braid(\g)$}

\bigskip

In \cite{WX}, it is proved that $\RC_\WX$ satisfies conditions ($\alpha$),
($\beta$) and ($\gamma$) of
Definition \ref{def:0.1}.
In \cite{GH2}, it is proved that it satisfies the first part of ($\delta$) of
this definition, namely $\RC_\WX$ induces the identity endomorphism
of $\m_{G^*\times G^*}/\m_{G^*\times G^*}^2$. Then $\RC_\WX -\Id$ induces a
map
$\m_{G^*\times G^*}/\m_{G^*\times G^*}^2 \to \m_{G^*\times
G^*}^2/\m_{G^*\times G^*}^3$, which we now compute.

\medskip

Identify $G^*$ with $\g^*$ using the exponential
map. We get, from (\ref{Req}), the  expansion at second order of the map $\RW$~:
\begin{align*}
\g^* \oplus \g^* &\to \g^* \oplus \g^*\\
(\xi,\eta)&\mapsto (\xi,\eta)+\left(
\ad^*(r_+(\eta))(\xi),\ad^*(r_-(\xi))(\eta)\right).
\end{align*}
View $(x,y) \in \g \oplus \g$ as a function of $\g^* \oplus \g^*\!$, taking
$(\xi,\eta)$ to
$\langle \xi,x\rangle + \langle \eta,y\rangle$. Then
$\RC_\WX (x,y)$ takes $(\xi,\eta)$ to
\begin{align*}
&\hskip-0.8cm\langle \xi+\ad^*(r_+(\eta))(\xi),x\rangle 
+\langle \eta+\ad^*(r_-(\xi))(\eta),y\rangle\\
=&\langle\xi,x\rangle+\langle\eta,y\rangle + \langle \xi,[r_+(\eta),x]\rangle
+\langle \eta,[r_-(\xi),y]\rangle\\
=&\langle\xi,x\rangle+\langle\eta,y\rangle +
\sum_i\langle a_i,\xi\rangle\langle\eta,[b_i,y]\rangle
+\sum_i\langle b_i,\eta\rangle\langle\xi,[a_i,x]\rangle\\
=&\langle\xi,x\rangle+\langle\eta,y\rangle +
\sum_i\langle \xi\otimes \eta,[a_i,x]\otimes b_i+a_i \otimes [b_i,y]\rangle +
\hbox {order } 3 \hbox{ in }
\!(\xi,\eta)
\end{align*}
where we set $r=\sum_ia_i \otimes b_i$, so that $r_+(\xi)=\sum_i\langle
b_i,\xi\rangle a_i$, and
$r_-(\xi)=\sum_i \langle a_i,\xi\rangle b_i$.
So
\begin{align*}
\left[\RC_\WX-\Id\right]~:~\g \oplus \g &\to \g \otimes \g\\
(x,y)&\mapsto [r,x\otimes1+1\otimes y].
\end{align*}
Then $\RC_\WX$ satisfies all the conditions of Definition \ref{def:0.1}.
\findem

\bigskip

\noindent {\bf - c - Unicity of braidings}

\bigskip

Let $\RC$ and $\RC'$ be two elements of $\Braid(\g)$.
We know that the maps
$[\RC-\Id]$ and $[\RC'-\Id]$~:
$\m_{G^*\times G^*}/\m_{G^*\times G^*}^2 \to 
\m_{G^*\times G^*}^2/\m_{G^*\times G^*}^3$ coincide, so
$(\RC-\RC')(\m_{G^*\times G^*})\subset
\m_{G^*\times G^*}^3$. Let us prove by induction over 
$k \geq 3$ that 
\begin{equation}
(\RC-\RC')(\m_{G^*\times G^*})\subset
\m_{G^*\times G^*}^k.
\label{id:k}
\end{equation}
As we have seen, (\ref{id:k}) holds for $k=3$.
Assume that it holds for some $k$ and let us prove it for 
$k+1$.
Let us set 
$\SC=\RC-\RC'$.
Then 
$\SC$ is a linear map
$\OC_{G^*\times G^*}\to m_{G^*\times G^*}^k$.
Moreover, we have for
$f,g \in \OC_{G^*\times G^*}$
\begin{equation}
\SC(fg)=\SC(f)\RC(g)+\RC(f)\SC(g)+\SC(f)\SC(g).
\label{quasi:der}
\end{equation}
Indentity (\ref{quasi:der}) allows to show by induction:
\begin{lemma}
\label{lemma:S}
For any $a \geq 1$, we have
$\SC(\m_{G^*\times G^*}^a)\subset\m _{G^*\times G^*}^{a+k-1}$.
\end{lemma}
\dem
This obviously holds when $a=1$.
\newline
Assume that we proved 
$\SC(\m_{G^*\times G^*}^a)\subset \m_{G^*\times G^*}^{a+k-1}$,
then for $f \in \m_{G^*\times G^*}^a$ and  $g\in \m_{G^*\times G^*}$,
\begin{multline*}
\SC(fg)=\SC(f)\RC(g)+\RC(f)\SC(g)+\SC(f)\SC(g)\\
\subset \!\m_{G^*\times G^*}^{a+k-1}\!\cdot\!
\m_{G^*\times G^*}+\m_{G^*\times G^*}^a\!\cdot\! \m_{G^*\times G^*}^k
+\m_{G^*\times G^*}^{a+k-1}\cdot\!\m_{G^*\times G^*}^k\!\subset\! \m_{G^*\times
G^*}^{a+k},
\end{multline*}
because
$a+2k-1 \geq a +k$. So
$\SC(\m_{G^*\times G^*}^{a+1}) \subset \m_{G^*\times G^*}^{a+k}$.
\findem

\smallskip
Let us now use the fact that $\OC_{G^*}$ is a topological Hopf algebra, equipped
with a decreasing Hopf filtration
$\OC_{G^*} \supset \m_{G^*} \supset \m_{G^*}^2\supset \cdots$.
The completion of the associated graded of $\OC_{G^*}$ is
a commutative and cocommutative Hopf algebra
$$\widehat{\gr}(\OC_{G^*})=\widehat{\oplus}_i \gr^i(\OC_{G^*})=\widehat{S}(\g).$$
$\OC_{G^*\times G^*}$ is also filtered and
$\widehat{\gr}(\OC_{G^*\times G^*})=\widehat{S}(\g)^{\bo 2}$.
Then Lemma \ref{lemma:S}, together with identity (\ref{quasi:der}), implies:
\begin{lemma}
Define $\gr(\SC)~:~
\widehat{\gr}\left(\OC_{G^*}^{\bo 2}\right)\to 
\widehat{\gr}\left(\OC_{G^*}^{\bo 2}\right)$ as the degree $k$
map such that
$\gr(\SC)$~:
$\m_{G^*\times G^*}^{a}/\m_{G^*\times G^*}^{a+1}\to
\m_{G^*\times G^*}^{a+k-1}/\m_{G^*\times G^*}^{a+k}$ is induced by $\SC$
for $a \geq 0$.
Then $\gr(\SC)$ is a derivation of degree
$k-1$ of 
$\gr(\OC_{G^*\times G^*})$.
\end{lemma}
\noindent
Comparing the analogues of the identities ($\gamma$) for $\RC$ and $\RC'$,
we get~:
\begin{align*}
\left(\SC^{1,3} \circ \RC^{2,3} +
\RC^{1,3} \circ \SC^{2,3}
+\SC^{1,3} \circ \SC^{2,3}\right)
\circ \left(\Delta_\OC \otimes \Id\right)&=\left(\Delta_\OC \otimes
\Id\right)\circ \SC,\\
\intertext{and }
\left(\SC^{1,3} \circ \RC^{1,2} +
\RC^{1,3} \circ \SC^{1,2}
+\SC^{1,3} \circ \SC^{1,2}\right)
\circ \left(\Id \otimes \Delta_\OC \right)&=
\left(\Id \otimes \Delta_\OC\right)\circ \SC.
\end{align*}
Both sides of each identity are algebra morphisms
$\OC_{G^*\times G^*}\to \OC_{G^* \times G^*\times G^*}$
taking 
$\m_{G^*\times G^*}^a$ to $\m_{G^* \times G^*\times G^*}^{a+k-1}$.
The associated graded morphisms are degree $k-1$ algebra morphisms
$\widehat{\gr}(\OC_{G^*\times G^*})\to \widehat{\gr}(\OC_{G^*\times G^*\times
G^*})$. The corresponding identities between these morphisms are
\begin{align}
\begin{split}
\label{id:gr}
\left(\gr\left(\SC\right)^{1,3} +\gr\left(\RC\right)^{2,3}\right)
\circ \left(\Delta_0 \otimes \Id\right)&=\left(\Delta_0 \otimes
\Id\right)\circ \gr(\SC),\\
\intertext{and }
\left(\gr\left(\SC\right)^{1,3}  +
\gr\left(\SC\right)^{1,2}\right)
\circ \left(\Id \otimes \Delta_0\right)&=
\left(\Id \otimes \Delta_0\right)\circ \gr(\SC),
\end{split}
\end{align}
where $\Delta_0$~: $\widehat{S}(\g)\to
\widehat{S}(\g) \bo \widehat{S}(\g)$ is the coproduct map
of $\gr(\SC)=\widehat{\gr}(\OC_{G^*}).$ These identities imply that the image of 
$\gr(\SC)$ is contained in $\Prim(\widehat{S}(\g))\otimes
\Prim(\widehat{S}(\g)).$
Since $\Prim(\widehat{S}(\g))=S^1(\g)=\gr^1(\OC_{G^*})$,
the image of $\gr(\SC)$ is therefore contained
in $\gr^1(\OC_{G^*})^{\otimes 2} \subset \gr^2(\OC_{G^*\times G^*})$.
Since the image of
$\gr( \SC)$ is also contained in
$\ho_{i\geq 3}\gr^i(\OC_{G^* \times G^*})$, we get 
$\gr(\SC)=0$.
It follows that $\SC(\m_{G^* \times G^*})\subset \m_{G^* \times G^*}^{k+1}$.
This proves the induction step of (\ref{id:k}).
Therefore
$\SC(\m_{G^* \times G^*})\subset \cap_{k\geq 0}\m_{G^* \times G^*}^k=0.$
Since $\SC$ is a derivation, we get
$\SC=0$.
Therefore 
$\RC=\RC'$.
This proves that $\Braid(\g)$ contains at most one element.
\hbox{}\findem
\vskip20pt

\centerline {\bf \S\; 3 \ Cohomological construction of $\rho$}

\stepcounter{section}\label{section:3}

\vskip20pt

Let $(\g,r)$ be a finite-dimensional quasitriangular Lie bialgebra. The purpose
of this section is to construct the unique element $\rho$ of $\Lift(\g)$ by
cohomological arguments, thus avoiding the use of associators.
Our main result is:
\begin{theorem}
\label{thm:existence}
$\Lift(\g)$ contains an element $\rho$.
\end{theorem}
\noindent
This result will be proved in Subsection 3.c.
In Subsection 3.a, we introduce variants and truncations of the sets $\Lift(\g)$
and $\Braid(\g)$.
Subsection 3.b contains the cohomological results allowing to construct
$\rho$ by successive approximations.

\bigskip

\noindent {\bf - a - Variants of the sets $\Braid(\g)$ and $\Lift(\g)$}

\bigskip

We denote by $\Braid'(\g)$ the set of all
Poisson automorphisms of
$\OC_{G^*\times G^*}$, satisfying conditions ($\alpha$), ($\gamma$) and
($\delta$) of Definition \ref{def:0.1}.
We denote by $\Lift'(\g)$ the set of all elements $\rho$ of 
$\OC_{G^*\times G^*}$, satisfying conditions 
($\alpha$), ($\gamma$) and
($\delta$) of Definition \ref{def:0.4}.
The map $\rho \mapsto \exp(V_\rho)$ then restrict to a map
$\Lift'(\g) \to \Braid'(\g)$. 

\smallskip

If $n$ is an integer, we define ${\Braid'\!}_{\leq n}\!(\g)$ (resp.,
${\Braid}_{\leq n}\!(\g)$) as the set of all Poisson automorphisms of
$\OC_{G^*\times G^*}/\m_{G^*\times G^*}^n$,
satisfying conditions ($\alpha$), ($\gamma$) and
($\delta$) (resp. ($\alpha$), ($\beta$), ($\gamma$) and
($\delta$)) of Definition \ref{def:0.1}, where $\OC_{(G^*)^k}$ is replaced
by $\OC_{(G^*)^k}/\m_{(G^*)^k}^n$, $k=1,2,3$.

\smallskip

Similarly, we define ${\Lift'\!}_{\leq n}\!(\g)$ (resp., $\Lift_{\leq n}\!(\g)$) 
as the set of all lifts 
$\rho \in \OC_{G^*\times G^*}/\m_{G^*\times G^*}^n$, satisfying conditions
($\alpha$), ($\gamma$) and
($\delta$) (resp. ($\alpha$), ($\beta$), ($\gamma$) and
($\delta$)) of Definition \ref{def:0.4},
where $\OC_{(G^*)^k}$ is replaced
by $\OC_{(G^*)^k}/\m_{(G^*)^k}^n$, $k=1,2,3$.
Then $\rho \mapsto \exp(V_\rho)$ defines a map
${\Lift'\!}_{\leq n}\!(\g)\to {\Braid'\!}_{\leq n}\!(\g)$.

\begin{lemma}
\label{lift:lift'}
We have:
\begin{enumerate}
\item The natural inclusions $\Lift(\g) \subset
\Lift'(\g)$, \ $\Braid(\g)\subset \Braid'(\g)$,\
${\Lift}_{\leq n}\!(\g)\subset {\Lift'\!}_{\leq n}\!(\g)$
and
${\Braid}_{\leq n}\!(\g)\subset {\Braid'\!}_{\leq n}\!(\g)$ are all equalities.
\item The set $\Braid_{\leq n}\!(\g)$ consists of only one element,
${\bar{\R}^{(n)}_{\WX}}$, which is the automorphism of
$\OC_{G^* \times G^*} /\m_{G^* \times G^*}^n$ induced
by the Weinstein-Xu automorphism.
\end{enumerate}
\end{lemma}
\dem
One can repeat the proof of the unicity part of Theorem \ref{thm:3}
to show that the sets ${\Braid'\!}_{\leq n}\!(\g)$,
${\Braid}_{\leq n}\!(\g)$ and ${\Braid'}(\g)$ all contain at most one element.
Since $\R_\WX$ is an element of $\Braid'( \g)$, we get $\Braid(\g)=
\Braid'(\g)=\{\R_\WX\}$.
In the same way,
the automorphism
${\bar{\R}^{(n)}_{\WX}}$ of $\OC_{G^* \times G^*} /\m_{G^* \times G^*}^n$ induced
by $\R_\WX$ is an element of
${\Braid'\!}_{\leq n}\!(\g)$ and
of ${\Braid}_{\leq n}\!(\g)$, so
${\Braid'\!}_{\leq n}\!(\g)={\Braid}_{\leq n}\!(\g)=\{\bar{\R}^{(n)}_{\WX}\}$.
This proves 2. and the equalities beween the sets of braidings of 1.

\smallskip

\noindent
Now $\Lift(\g)$ is defined as the preimage of $\Braid(\g)$ by
the map
\begin{align*}
\exp~:\Lift'(\g) &\to \Braid'(\g)\\
\rho & \mapsto \exp(V_\rho);
\end{align*}
similarly, ${\Lift\!}_{\leq n}\!(\g)$ is the preimage
of
${\Braid\!}_{\leq n}\!(\g)$ by the map $\exp$~: 
${\Lift'\!}_{\leq n}\!(\g)$ $ \to {\Braid'\!}_{\leq n}(\g)$.
So we get $\Lift(\g)=\Lift'(\g)$ and 
${\Lift\!}_{\leq n}(\g)={\Lift'\!}_{\leq n}(\g)$.
\findem

\bigskip

\noindent {\bf - b - A map ${\Lift'\!}_{\leq n}\!(\g) \to {\Lift'\!}_{\leq
n+1}\!(\g)$}

\bigskip

We have canonical projection maps
${\Lift'\!}_{\leq n}\!(\g)\xrightarrow{\pi_{n-1}} {\Lift'\!}_{\leq n-1}\!(\g)
\to \cdots$. Then
$$\Lift'(\g)=\lim_{\Tto_n}({\Lift'\!}_{\leq n}\!(\g)).$$
To construct an element of $\Lift'(\g)$, we will therefore
construct a sequence of maps
$$\lambda_n~:~{\Lift'\!}_{\leq n}\!(\g)\to {\Lift'\!}_{\leq n+1}\!(\g),~n \geq 3,$$
such that $\pi_n \circ \lambda_n=\Id$.

\smallskip

\noindent Let $\rho_n\in \OC_{G^* \times G^*} /\m_{G^* \times G^*}^n$
be an element of ${\Lift'\!}_{\leq n}\!(\g)$. We have then
$$(\epsilon \otimes \Id)(\rho_n)=(\Id \otimes \epsilon)(\rho_n)=0,$$
$$(\Delta \otimes \Id)(\rho_n)=\rho_n^{1,3} \star \rho_n^{2,3},~
(\Id \otimes \Delta)(\rho_n)=\rho_n^{1,3} \star \rho_n^{1,2},$$
$$[\rho_n]=r.$$
Let us take a lift $\widetilde{\rho}_n
\in \OC_{G^* \times G^*}/
\m_{G^*\times G^*}^{n+1}$ of $\rho_n$ such that
$(\epsilon \otimes \Id)(\widetilde{\rho}_n)=
(\Id \otimes \epsilon)(\widetilde{\rho}_n)=0.$
Set
\begin{align}
\begin{split}
\label{alpha:beta}
\alpha&=(\Delta \otimes \Id)(\widetilde\rho_n)-
\widetilde\rho_n^{1,3} \star \widetilde{\rho}_n^{2,3},\\
\beta&=(\Id \otimes \Delta)(\widetilde\rho_n)=
\widetilde\rho_n^{1,3} \star \widetilde{\rho}_n^{1,2}.
\end{split}
\end{align}
Then $\alpha,\beta \in \m_{(\OC^*)^3}^n/
\m_{(\OC^*)^3}^{n+1}$. Moreover
\begin{align*}
(\epsilon \otimes \Id \otimes \Id)(\alpha)&=(\Id \otimes \,\epsilon \otimes
\Id)(\alpha)=(\Id \otimes \Id \otimes \, \epsilon)(\alpha)=0\\
\hbox{and~}
(\epsilon \otimes \Id \otimes \Id)(\beta)&=(\Id \otimes\, \epsilon \otimes
\Id)(\beta)=(\Id \otimes \Id \otimes \, \epsilon)(\beta)=0.
\end{align*}
Let $\sigma$ be an element of
$\m_{G^* \times G^*}^n/\m_{G^* \times G^*}^{n+1}$.
Set $\rho_{n+1}=\widetilde{\rho}_n+\sigma$.
Then $\rho_{n+1}$ belongs to 
${\Lift'\!}_{\leq n+1}\!(\g)$ if and only if:
\begin{equation}
\label{counit}
(\epsilon \otimes \Id)(\sigma)=
(\Id \otimes \,\epsilon)(\sigma)=0,
\end{equation}
\begin{equation}
\label{coboundary}
(\dif \otimes \Id)(\sigma)=-\alpha,~
(\Id \otimes \dif)(\sigma)=-\beta.
\end{equation}

Here, we identify $\m_{G^* \times G^*}^n/\m_{G^*
\times G^*}^{n+1}$ with $S^n(\g \oplus \g)$ and 
$\m_{(\OC^*)^3}^n/
\m_{(\OC^*)^3}^{n+1}$ with $S^n(\g \oplus \g \oplus \g)$.
Then the map $\dif$~: 
$S^\cdot(\g) \to S^\cdot(\g \oplus \g)$ is the Hochschild coboundary map,
taking $f$ to $\Delta_0(f)-f\otimes 1 - 1 \otimes f$ ($\Delta_0$
is the cocommutative coproduct of $S^\cdot(\g)$).
Identities (\ref{counit}) and (\ref{coboundary}) follow from the identities 
\begin{equation}
\label{f:g:h}
f\star (h+g)=f\star h+g,\hskip0.6cm (f+g)\star h= f\star h +g,
\end{equation}
when $f,h \in \m_{(G^*)^k}^2/\m_{(G^*)^k}^{n+1}$
and $g \in \m_{(G^*)^k}^n/\m_{(G^*)^k}^{n+1}$.

\medskip

\noindent
Let us now recall some results of co-Hochschild cohomology.
Let $\dif^{(2)}$~:
$S^\cdot(\g \oplus \g)\to 
S^\cdot(\g \oplus \g \oplus \g)$ be defined by
$\dif^{(2)}(f)=(\dif\otimes \Id)(f) -(\Id \otimes \dif)(f)$ (we identify
$S^\cdot(\g^{\oplus k})$ with $S^\cdot(\g)^{\otimes k}$).
Then $\dif^{(2)} \circ \dif=0$.
The cohomology group $H^1_\cohoch =
\Ker(\dif^{(2)})/\Im(\dif)$  identifies with $\wedge^2(\g)$.
The canonical map $\Ker(\dif^{(2)}) \to \wedge^2(\g)$ is given by the
antisymmetrization $f \mapsto f - f^{2,1}$.
The $0$-th cohomology group 
$H^0_\cohoch=\Ker(\dif)$ is equal to $\g$. We then prove:
\begin{lemma}
\label{lemma:sigma}
There exists a solution $\sigma \in \m_{G^* \times G^*}^n/\m_{G^* \times
G^*}^{n+1}$ of equations 
(\ref{counit}) and (\ref{coboundary}) if and only if $\alpha,\beta$
satisfy the equations:
\begin{align}
\label{conds:1}
(\Id \otimes \Id \otimes \dif)(\alpha)&=(\dif \otimes \Id \otimes \Id)(\beta),\\
\label{conds:2}
(\dif^{(2)}\otimes \Id)(\alpha)&=(\Id \otimes \dif^{(2)}(\beta)=0,\\
\label{conds:3}
\alpha^{1,2,3}=\alpha^{2,1,3}, \; & \; 
\beta^{1,2,3}=\beta^{1,3,2}.\end{align}
If these conditions are satisfied, then the solution is unique.
\end{lemma}
\dem
Assume that $\sigma$ exists.
Then both sides of
(\ref{conds:1})
are equal to $-(\dif \otimes \dif)(\sigma)$,
so we have (\ref{conds:1}).
(\ref{conds:2}) follows from $\dif^{(2)} \circ \dif=0$
and (\ref{conds:3}) follows from the fact that the image of
$\dif$~: $S^\cdot(\g) \to S^\cdot(\g \oplus \g)$ is contained in the subspace of
invariants under the permutation of both summands of $\g \oplus \g$.
So (\ref{conds:1}), (\ref{conds:2}) and (\ref{conds:3}) are satisfied.

\smallskip

\noindent Assume now that these identities are satisfied.
The equalities $(\dif^{(2)} \otimes \Id)(\alpha)=0$ and
$\alpha^{1,2,3}=\alpha^{2,1,3}$ imply that there exists
$\sigma' \in S^n(\g \oplus \g)$, such
that 
$(\dif \otimes \Id)(\sigma')=-\alpha$.
In the same way, there exists $\sigma'' \in S^n(\g \oplus \g)$, such
that 
$(\Id \otimes \dif)(\sigma'')=-\beta$.
$\sigma'$ is well-defined only up to addition of an element of
$\g \otimes S^{n-1}(\g)$, and $\sigma''$ is well-defined up to addition
of an element of $S^{n-1}(\g) \otimes \g$.
Now (\ref{conds:1}) implies that
$(\dif \otimes \dif)(\sigma'-\sigma'')=0$.
Since $\Ker(\dif)=\g$, we get
$\sigma'-\sigma'' \in \g \otimes
S^{n-1}(\g) + S^{n-1}(\g) \otimes \g$.
Let $\sigma'-\sigma''={\sigma'}_0 + {\sigma''}_0$,
${\sigma'}_0\in \g \otimes S^{n-1}(\g)$,
${\sigma''}_0\in  S^{n-1}(\g) \otimes \g$.
Set
$\sigma=\sigma'-{\sigma'}_0=\sigma''+{\sigma''}_0$.
Then
\begin{align*}
(\dif \otimes \Id)(\sigma)&=(\dif \otimes \Id)(\sigma')=-\alpha\\
\hbox{and~}~(\Id \otimes \dif)(\sigma)&=(\Id \otimes \dif)(\sigma'')=-\beta.
\end{align*}
So equations (\ref{conds:1}), (\ref{conds:2}) and (\ref{conds:3}) imply
the existence of a
solution $\sigma$ of (\ref{coboundary}).
Then $(\epsilon \otimes \, \epsilon)\circ \dif = -\epsilon$, so
$$(\epsilon \otimes \Id)(\sigma)=-(\epsilon \otimes \, \epsilon)
\circ (\dif \otimes \Id)(\sigma)=(\epsilon \otimes \, \epsilon
\otimes \Id)(\alpha)=0,$$
and in the same way $(\Id \otimes \, \epsilon)(\sigma)=0$.
So $\sigma$ satisfies also
equation (\ref{counit}).

\smallskip

\noindent The unicity of $\sigma$ then follows from the fact that
$\Ker(\dif \otimes \Id) \cap \Ker(\Id \otimes \dif) =\g \otimes \g \subset
S^2(\g \oplus \g)$, so the intersection of
$\Ker(\dif \otimes \Id)\cap \Ker(\Id \otimes \dif)$ with $S^n(\g \oplus \g)$
is zero as $n\geq 3$.
\findem

\begin{proposition}
\label{ids:OK}
The elements $\alpha$ and $\beta$ defined by (\ref{alpha:beta}) satisfy the
identities (\ref{conds:1}), (\ref{conds:2}) and (\ref{conds:3}).
\end{proposition}
\dem
Apply $\Delta \otimes \Id \otimes \Id -\Id \otimes \Delta \otimes \Id$ to
the identity
\begin{equation}
\label{id:1}
(\Delta \otimes \Id)(\widetilde{\rho}_n)=\widetilde{\rho}_n^{1,3} \star
\widetilde{\rho}_n^{2,3} + \alpha.
\end{equation}
This yields an identity in $\OC_{(G^*)^4}/\m_{(G^*)^4}^{n+1}$.
Its left side vanishes since
$(\Delta \otimes \Id -\Id \otimes \Delta) \circ \Delta=0$.
Using again (\ref{id:1}),
we get
$$0=(\widetilde{\rho}_n^{1,4} \star \widetilde{\rho}_n^{2,4} + \alpha^{1,2,4})\star
\widetilde{\rho}_n^{3,4}+\alpha^{12,3,4} - \widetilde{\rho}_n^{1,4} \star (\widetilde{\rho}_n^{2,4} \star
\widetilde{\rho}_n^{3,4} + \alpha^{2,3,4})-\alpha^{1,23,4}.$$
Now identities (\ref{f:g:h}) yield:
$$
0=
(\widetilde{\rho}_n^{1,4} \star \widetilde{\rho}_n^{2,4}\star 
\widetilde{\rho}_n^{3,4}
+ \alpha^{1,2,4})+\alpha^{12,3,4}\!- (\widetilde{\rho}_n^{1,4} 
\star \widetilde{\rho}_n^{2,4} \star
\widetilde{\rho}_n^{3,4} \!+ \alpha^{2,3,4})-\alpha^{1,23,4}
$$
that is 
$$
(\dif^{(2)}\otimes \Id)(\alpha)=0.
$$
Applying $\Id \otimes \Delta \otimes \Id - \Id \otimes \Id \otimes \Delta$ to
the identity
\begin{equation}
\label{id:2}
(\Id \otimes \Delta)(\widetilde{\rho}_n)=\widetilde{\rho}_n^{1,3} \star
\widetilde{\rho}_n^{1,2}+\beta,
\end{equation}
we get in the same way
$$(\Id \otimes \otimes \dif^{(2)})(\beta)=0.$$
So $\alpha$ and $\beta$ satisfy (\ref{conds:2}).

\bigskip

Apply now $\Id \otimes \Id \otimes \Delta$ to (\ref{id:1}),
$\Delta \otimes \Id \otimes \Id$ to (\ref{id:2}), and
substract the resulting equalities.
Using again (\ref{id:1}) and (\ref{id:2}), we get
\begin{align*}
0=&(\widetilde{\rho}_n^{1,4} \star \widetilde{\rho}_n^{1,3}
+ \beta^{1,3,4})\star (\widetilde{\rho}_n^{2,4} \star \widetilde{\rho}_n^{2,3}
+ \beta^{2,3,4})+\alpha^{12,3,4}\\
&-
(\widetilde{\rho}_n^{1,4} \star \widetilde{\rho}_n^{2,4}+\alpha^{1,2,4})
\star(\widetilde{\rho}_n^{1,3}
\star \widetilde{\rho}_n^{2,3}+\alpha^{1,2,3})-\beta^{1,2,34}.
\end{align*}
Using again identities (\ref{f:g:h}), and the fact that $\widetilde{\rho}_n^{1,3}
\star \widetilde{\rho}_n^{2,4}=\widetilde{\rho}_n^{2,4}\star \widetilde{\rho}_n^{1,3}$,
we get
$$\alpha^{12,3,4}-\alpha^{1,2,3}-\alpha^{1,2,4}=\beta^{1,2,34}-\beta^{1,2,3}
-\beta^{1,2,4},$$
that is $(\dif \otimes \Id \otimes \Id)(\alpha)=(\Id \otimes \Id \otimes
\dif)(\beta)$. So $\alpha$ and $\beta$ satisfy (\ref{conds:1}).

\bigskip

To prove that they also satisfy (\ref{conds:3}), let us set
$$ 
\psi =
\widetilde{\rho}_n^{1,2}\star \widetilde{\rho}_n^{1,3}\star 
\widetilde{\rho}_n^{2,3}
-\widetilde{\rho}_n^{2,3}\star \widetilde{\rho}_n^{1,3} \star \widetilde{\rho}_n^{1,2}
$$
and let us prove:
\begin{lemma}
\label{Y:0}
We have $\psi=0$.
\end{lemma}
\noindent{\sc Proof of Lemma.}
Since $\rho_n^{1,2}\star {\rho}_n^{1,3}\star {\rho}_n^{2,3}=
\rho_n^{1,2}\star \rho_n^{12,3}=
\rho_n^{21,3}\star \rho_n^{1,2}=
{\rho}_n^{2,3}\star {\rho}_n^{1,3} \star {\rho}_n^{1,2}$,
the class of $\psi$ in 
$\OC_{(G^*)^3}/\m_{(G^*)^3}^n$ is zero, so $\psi \in
\m_{(G^*)^3}^n/\m_{(G^*)^3}^{n+1}$.
We identify $\psi$ with an element of $S^n(\g \oplus \g \oplus \g)$.
Then
\begin{align*}
\psi^{12,3,4} =& \widetilde{\rho}_n^{12,3}\star \widetilde{\rho}_n^{12,4}\star
\widetilde{\rho}_n^{3,4} - \widetilde{\rho}_n^{3,4} \star 
\widetilde{\rho}_n^{12,4} \star \widetilde{\rho}_n^{12,3}
\\
= &(\widetilde{\rho}_n^{1,3} \star \widetilde{\rho}_n^{2,3}
+\alpha^{1,2,3})\star
(\widetilde{\rho}_n^{1,4} \star \widetilde{\rho}_n^{2,4}+\alpha^{1,2,4})\star
\widetilde{\rho}_n^{3,4}\\
&-\widetilde{\rho}_n^{3,4}\star (\widetilde{\rho}_n^{1,4}
 \star \widetilde{\rho}_n^{2,4}+\alpha^{1,2,4})\star
 (\widetilde{\rho}_n^{1,3} \star \widetilde{\rho}_n^{2,3}+\alpha^{1,2,3})
\\  =& \widetilde{\rho}_n^{1,3}\star \widetilde{\rho}_n^{2,3}
\star \widetilde{\rho}_n^{1,4}
 \star \widetilde{\rho}_n^{2,4}\star \widetilde{\rho}_n^{3,4}
 -\widetilde{\rho}_n^{3,4}\star \widetilde{\rho}_n^{1,4}\star 
\widetilde{\rho}_n^{2,4}
 \star \widetilde{\rho}_n^{1,3}\star \widetilde{\rho}_n^{2,3}
\\
&\hbox{~}\hskip7.6cm \hbox{by virtue of (\ref{f:g:h})}
\\
=&\widetilde{\rho}_n^{1,3}\star \widetilde{\rho}_n^{1,4}\star 
( \widetilde{\rho}_n^{3,4} \star \widetilde{\rho}_n^{2,4}\star 
\widetilde{\rho}_n^{2,3} + \psi^{2,3,4})
\\
& + ( - \widetilde{\rho}_n^{1,3} \star \widetilde{\rho}_n^{1,4}
\star \widetilde{\rho}_n^{3,4} + \psi^{1,3,4}) \star
 \widetilde{\rho}_n^{2,4}\star \widetilde{\rho}_n^{2,3}
\\
&\hbox{~}\hskip1.6cm \hbox{since 
$\widetilde\rho^{1,3} \star \widetilde\rho^{2,4} 
= \widetilde\rho^{2,4} \star \widetilde\rho^{1,3}$ 
and 
$\widetilde\rho^{1,4} \star \widetilde\rho^{2,3} 
= \widetilde\rho^{2,3} \star \widetilde\rho^{1,4}$ }
\\
 =&\psi^{1,3,4}+\psi^{2,3,4}\hskip5.4cm \hbox{by virtue of (\ref{f:g:h})}.
\end{align*}
We have also
\begin{align*}
\psi^{1,23,4}=&(\widetilde{\rho}_n^{1,3} \star \widetilde{\rho}_n^{1,2}
+\beta^{1,2,3})\star
\widetilde{\rho}_n^{1,4}\star
(\widetilde{\rho}_n^{2,4} \star \widetilde{\rho}_n^{3,4}+\alpha^{2,3,4})\\
&-(\widetilde{\rho}_n^{2,4}
 \star \widetilde{\rho}_n^{3,4}+\alpha^{2,3,4})\star
 \widetilde{\rho}_n^{1,4} \star 
 (\widetilde{\rho}_n^{1,3} \star \widetilde{\rho}_n^{1,2}+\beta^{1,2,3})\\
 =&\widetilde{\rho}_n^{1,3}\star \widetilde{\rho}_n^{1,2}\star 
\widetilde{\rho}_n^{1,4}
 \star \widetilde{\rho}_n^{2,4}\star \widetilde{\rho}_n^{3,4}
 -\widetilde{\rho}_n^{2,4}\star \widetilde{\rho}_n^{3,4}\star 
\widetilde{\rho}_n^{1,4}
 \star \widetilde{\rho}_n^{1,3}\star \widetilde{\rho}_n^{1,2}\\
 =&\widetilde{\rho}_n^{1,3}\star 
(\widetilde{\rho}_n^{2,4} \star 
\widetilde{\rho}_n^{1,4}\star \widetilde{\rho}_n^{1,2} + \psi^{1,2,4})
 \star \widetilde{\rho}_n^{3,4}\\
 &+\widetilde{\rho}_n^{2,4}\star 
( - \widetilde{\rho}_n^{1,3}
 \star \widetilde{\rho}_n^{1,4}\star \widetilde{\rho}_n^{3,4}
+ \psi^{1,3,4}) 
\star  \widetilde{\rho}_n^{1,2}\\\
 =&\psi^{1,2,4}+\psi^{2,3,4}.
\end{align*}
In the same way, one proves that
$\psi^{1,2,34}=\psi^{1,2,3}+\psi^{1,2,4}$. Therefore, we get
$$(\dif \otimes \Id \otimes \Id)(\psi)=
(\Id \otimes \dif \otimes \Id)(\psi)=(\Id \otimes \Id \otimes \dif)(\psi)=0,$$
so $\psi \in \g \otimes \g \otimes \g$. When $n >3$, this implies $\psi=0$.
When $n=3$, $\psi$ is equal
to $[r^{1,2},r^{1,3}]+[r^{1,2},r^{2,3}]+[r^{1,3},r^{2,3}]$ and is also zero.
\findem

\medskip

\noindent{\sc End of proof of Proposition \ref{ids:OK}}
Let us now prove that $\alpha$ and $\beta$ satisfy equation (\ref{conds:3}).
For this, we first prove: 

\begin{lemma}
We have 
\begin{equation} \label{QT:id}
\widetilde \rho_n^{1,2} \star \widetilde \rho_n^{12,3}
= \widetilde\rho_n^{21,3} \star \widetilde \rho_n^{1,2} 
\end{equation}
(equality in $\OC_{(G^*)^3}/\m_{(G^*)^3}^{n+1}$). 
\end{lemma}

\noindent{\sc Proof of Lemma.} We will prove that if 
$$
f \in ( \m_{G^*} \bar\otimes \m_{G^*}) / 
\big( \m_{G^*\times G^*}^{n+1}\cap ( \m_{G^*} \bar\otimes \m_{G^*})  \big), 
$$
the equality 
\begin{equation} \label{eq:rho}
\widetilde{\rho}+n^{1,2} \star f^{12,3} 
= f^{21,3} \star \widetilde{\rho}_n^{1,2}
\end{equation}
holds in $\OC_{(G^*)^3}/\m_{(G^*)^3}^{n+1}$. 

There exist element $f'_i,f''_i$ of $\m_{G^*}$, such that $f$ is equal 
to the class of $\sum_i f'_i \otimes f''_i$. Then we have 
$\widetilde\rho_n^{1,2} \star (f'_i)^{12} = (f'_i)^{21} \star 
\widetilde\rho_n^{1,2}$ (equality in $\OC_{(G^*)^3}/\m_{(G^*)^3}^{n}$), 
because $\rho_n \in \operatorname{Lift}'_{\leq n}(\g) = 
\operatorname{Lift}'_{\leq n}(\g)$ by virtue of Lemma 
\ref{lift:lift'}. Tensoring this identity with $f''_i\in \m_{G^*}$, 
we get 
$$
\big( \widetilde\rho_n^{1,2} \star (f'_i)^{12}\big) (f''_i)^{3}
= 
\big( (f'_i)^{21} \star \widetilde\rho_n^{1,2} \big) (f''_i)^{3} ,
$$
an equality in 
$$
\big( \OC_{G^*\times G^*} \bar\otimes \m_{G^*}\big) 
/ \big( \m_{G^*\times G^*}^n \bar\otimes \m_{G^*}\big)
$$ 
and therefore also in 
$\OC_{(G^*)^3} / \m_{(G^*)^3}^{n+1}$.  Summing over all indices $i$, 
we get (\ref{eq:rho}), which implies (\ref{QT:id}) by taking 
$f = \widetilde{\rho}_n$. \hfill \qed \medskip 

Plugging (\ref{id:1}) into (\ref{QT:id}), we get
$$
\widetilde{\rho}_n^{1,2} \star (\widetilde{\rho}_n^{1,3} \star 
\widetilde{\rho}_n^{2,3}
+ \alpha^{1,2,3})=(\widetilde{\rho}_n^{2,3} \star \widetilde{\rho}_n^{1,3}
+ \alpha^{2,1,3})
\star \widetilde{\rho}_n^{1,2}.
$$
Then (\ref{f:g:h}) yields:
$$\alpha^{2,1,3}-\alpha^{1,2,3}=\widetilde{\rho}_n^{1,2}
\star \widetilde{\rho}_n^{1,3}\star \widetilde{\rho}_n^{2,3}
-  \widetilde{\rho}_n^{2,3}
\star \widetilde{\rho}_n^{1,3}
\star \widetilde{\rho}_n^{1,2},$$
so Lemma \ref{Y:0} gives:
$$\alpha^{2,1,3}-\alpha^{1,2,3}=0.$$
In the same way, we have
$$\widetilde{\rho}_n^{2,3} \star \widetilde{\rho}_n^{1,23}=\widetilde{\rho}_n^{1,32}
\star \widetilde{\rho}_n^{2,3},$$
so by (\ref{id:2}), we get
$$\widetilde{\rho}_n^{2,3} \star (\widetilde{\rho}_n^{1,3} \star \widetilde{\rho}_n^{1,2}
+ \beta^{1,2,3})=(\widetilde{\rho}_n^{1,2} \star \widetilde{\rho}_n^{1,3}
+ \beta^{1,3,2})
\star \widetilde{\rho}_n^{2,3},$$
so $\beta^{1,2,3}-\beta^{1,3,2}=\widetilde{\rho}_n^{1,2}
\star \widetilde{\rho}_n^{1,3}\star \widetilde{\rho}_n^{2,3}
-  \widetilde{\rho}_n^{2,3}
\star \widetilde{\rho}_n^{1,3}
\star \widetilde{\rho}_n^{1,2},$
so by Lemma \ref{Y:0}, $\beta^{1,2,3}-\beta^{1,3,2}=0$.
So $\alpha$ and $\beta$ satisfy equation (\ref{conds:3}).
\findem

\medskip

\noindent Let us now construct the map $\lambda_n$~:
${\Lift'\!}_{\leq n}\!(\g) \to {\Lift'\!}_{\leq
n+1}\!(\g)
$. Let $\rho \mapsto \widetilde{\rho}$ be any map
\begin{multline*}
\{\rho \in \OC_{G^* \times G^*}/\m_{G^* \times G^*}^n|
(\epsilon \otimes \Id)(\rho)=(\Id \otimes \, \epsilon)(\rho)=0\}\\
\to \{\widetilde{\rho} \in \OC_{G^* \times G^*}/\m_{G^* \times G^*}^{n+1}|
(\epsilon \otimes \Id)(\widetilde{\rho})=(\Id \otimes \, \epsilon)
(\widetilde{\rho})=0\},
\end{multline*}
which is a section of the canonical projection map (we may take $\rho \mapsto
\widetilde{\rho}$ linear).
If $\rho_n \in {\Lift'\!}_{\leq n}\!(\g)$,
we define $\alpha$ and $\beta$ by (\ref{alpha:beta}).
Then Proposition
\ref{ids:OK} and Lemma \ref{lemma:sigma} allow us to construct a unique element
$\sigma \in \m_{G^* \times G^*}^n/\m_{G^* \times G^*}^{n+1}$,
such that $\widetilde{\rho}_n+\sigma \in {\Lift'\!}_{\leq n+1}\!(\g)$.
We then set
$$\lambda_n(\rho_n)=\widetilde{\rho}_n+\sigma.$$
This defines the desired map
$\lambda_n$~:
${\Lift'\!}_{\leq n}\!(\g) \to {\Lift'\!}_{\leq
n+1}\!(\g)
$.

\bigskip

\noindent {\bf - c - {\sc Proof of Theorem \ref{thm:existence}}}

\bigskip

$r \in \g \otimes \g$ defines an element 
$\rho_3$ of ${\Lift'\!}_{\leq 3}\!(\g)$.
Applying to it 
$\lambda_3$, $\lambda_4,\dots,$ we define a sequence $\rho_n$ of
elements of ${\Lift'\!}_{\leq n}\!(\g)$, and
therefore an element
$\rho \in \Lift'(\g)$.
According to Lemma \ref{lift:lift'}, $\rho$ is then an element
of $\Lift(\g)$.
\findem

\vskip20pt

\centerline {\bf \S\; 4 \ Construction of universal lifts (proof of Theorem
\ref{thm:1.10})}

\stepcounter{section}\label{section:4}

\vskip20pt

Theorem \ref{thm:1.10} can be proved in the same way as its ``non-universal''
counterpart Theorem \ref{thm:0.8}:
\begin{enumerate}
\item the unicity part is proved using the same argument;
\item the existence part can be proved either using the map $\Quant(\KM)
\to \Lift_\univ$, and the nonemptiness of $\Quant(\KM)$ (see \cite{EK});
or it can be proved following the arguments of Section 4.
\end{enumerate}

\vskip20pt

\centerline {\bf \S\hskip0.1cm 5 Appendix: a commutative diagram related to QFSH
algebras}

\stepcounter{section}\label{section:5}

\vskip20pt

The aim of this section is to prove the following lemma:
\begin{lemma}
Let $\sigma$  be an arbitrary element of $\m_\hbar$
and $[\sigma]$ be its class in
$\m_\hbar/\hbar \m_\hbar + \m_\hbar^2$.
Since $\m_\hbar/\hbar \m_\hbar + \m_\hbar^2$ identifies
with $\g$,
$[\sigma]\in \g$.
Since $\m_\hbar \subset \hbar U_\hbar(\g)$,
$\sigma$ is an element of  $\hbar U_\hbar(\g)$.
Then $\left({{\sigma}\over{\hbar}} \mod \hbar\right)$ is an element of
$U(\g)$.
We have the following identity in $U(\g)$~:
$$\left({{\sigma}\over{\hbar}}\mod \hbar \right)=[\sigma].$$
\label{gav2}
\end{lemma}
\noindent This lemma clearly implies Lemma \ref{gav}: if 
$\sigma=\sum_i \sigma_i^1 \otimes
\sigma_i^2 \in \m_\hbar \bo \m_\hbar$ 
satisfies the hypothesis of Lemma \ref{gav},
then $[\sigma]=\sum_i [\sigma_i^1] \otimes
[\sigma_i^2] \in \g^{\otimes 2}$ and
\begin{align*}
\left.{{\sigma}\over{\hbar^2}}\right|_{U_{\hbar}^{\ho 2}\to U^{\otimes 2}}&=
\sum_i \left.\left( {{\sigma_i^1}\over{\hbar}} \otimes
{{\sigma_i^2}\over{\hbar}}\right)\right|_{U_{\hbar}^{\ho 2}\to U^{\otimes 2}}\\
&=
\sum_i \left.{{\sigma_i^1}\over{\hbar^2}}\right|_{U_{\hbar}\to U}
\otimes \left.{{\sigma_i^2}\over{\hbar^2}}\right|_{U_{\hbar}\to U}\\
&=\sum_i [\sigma_i^1] \otimes
[\sigma_i^2] &\text{(by Lemma \ref{gav2})}\\
&=[\sigma].
\end{align*}
We use the notation $\left.x\right|_{U_{\hbar}^{\ho k}\to
U^{\otimes k}}$ for $(x \mod \hbar)$, when $x \in U_{\hbar}(\g)^{\ho k}$.
\medskip

\noindent More generally, in this section, we will consider a Lie bialgebra $(\g,\delta)$
over a field $\KM$,
$(U_\hbar(\g),\Delta)$ a quantization of $U(\g)$ and the subalgebra
$U_\hbar(\g)' \subset U_\hbar(\g)$, where
$$U_\hbar(\g)'=\{x\in U_\hbar(\g)|~\forall n \in \NM,~ \delta^{(n)}(x) 
\in \hbar^n 
U_\hbar(\g)^{\widehat{\otimes} n}\}$$ 
(where $\delta^{(n)}=\left(\Id - \eta \circ \varepsilon\right)^{\otimes
n}\circ \Delta^{(n)}$). 
The definition of $U_\hbar(\g)'$ yields, when $n=1$ or $2$: 
\begin{lemma}
\label{lemma:reduc}
Let $f$ be an element of $U_\hbar(\g)'$. We have:
\begin{itemize}
\item $f -\varepsilon(f) \in \hbar U_\hbar(\g)$
\item $\left.{{f-\varepsilon(f)}\over{\hbar}}\right|_{U_{\hbar}\to U}\in
\g \subset U(\g)$
\end{itemize}
\end{lemma}
\noindent
According to a theorem of Drinfeld (see \cite{Dr} and also
\cite{Ga}) the quantized formal series Hopf algebra $U_\hbar(\g)'$ is a
quantization of the function algebra $\OC_{G^*}$.
The projection
$U_\hbar(\g)'/\hbar U_\hbar(\g)'\to \OC_{G^*}\simeq U(\g^*)^*$
may be described as follows:
\begin{theorem}
\label{theo:lambda}
For $f \in U_\hbar(\g)'$, let
$\left.f\right|_{\OC_{\hbar}\to
\OC}$ be its class in
$U_\hbar(\g)'/\hbar U_\hbar(\g)'$.
There exists a unique Hopf pairing
$U(\g^*)\otimes \left(U_\hbar(\g)'/\hbar U_\hbar(\g)'\right)$ $\to \KM$,
such that 
$$\forall \xi \in \g^*, \forall f \in U_\hbar(\g)'~
\left\langle \xi ,\left.f\right|_{\OC_{\hbar}\to
\OC}\right\rangle=\left\langle  \xi,
\left.{{f-\varepsilon(f)}\over{\hbar}}\right|_{U_{\hbar}\to U}\right\rangle_{\g\times
\g^*}.$$
This pairing induces an isomorphism
$U_\hbar(\g)'/\hbar U_\hbar(\g)'\To_\lambda^\sim U(\g^*)^*$.
\end{theorem}
\noindent
We can now reformulate Lemma \ref{gav2} (in this section, we will use the
notation $(H)_0$ for the maximal ideal $\Ker(\varepsilon)$ of a Hopf algebra
$H$):
\begin{proposition}
The following diagram commutes
$$\begin{matrix}
{\left(U_\hbar(\g)'\right)_0/\left(\hbar
\left(U_\hbar(\g)'\right)_0+\big(U_\hbar(\g)'\right)_0^2\big{)}}
&\xrightarrow{\Simto^{\scriptstyle{\rm{(a)}}}}&\left(U(\g^*)^*\right)_0/
\left(U(\g^*)^*\right)_0^2\cr
&&\cr
{\scriptstyle{\rm{(b)}}}\downarrow\hskip-2cm&&\downarrow{\scriptstyle{\rm{(d)}}}\cr
{\hbar U_\hbar(\g)/\hbar^2 U_\hbar(\g)\hskip-3cm}
&\xrightarrow{\Simto^{\scriptstyle{\rm{(c)}}}}
&{U(\g)\Hoo_{\cang}\g\hskip2cm}
\end{matrix}$$
where $\cange$~: $\g \hookrightarrow U(\g)$ is the canonical injection and the other maps are given by:
\begin{description}
\item[(a)] is the composed map
\begin{multline*}
{\left(U_\hbar(\g)'\right)_0/\big(\hbar
\left(U_\hbar(\g)'\right)_0+\left(U_\hbar(\g)'\right)_0^2\big)}\\
{\To^\sim
\{\{\left(U_\hbar(\g)'\right)_0/\hbar
\left(U_\hbar(\g)'\right)_0\}\}/
\{\{\left(U_\hbar(\g)'\right)_0/\hbar
\left(U_\hbar(\g)'\right)_0\}\}^2}\\
\to \left(U(\g^*)^*\right)_0/
\left(U(\g^*)^*\right)_0^2
\end{multline*} 
where the last map is induced by
$\lambda$~: $U_\hbar(\g)'/\hbar U_\hbar(\g)'\To_\lambda^\sim U(\g^*)^*$
(see Theorem \ref{theo:lambda}),
\item[(b)] is the quotient map of the injection~:
$\left(U_\hbar(\g)'\right)_0 \hookrightarrow \hbar U_\hbar(\g)$
with respect to the ideals $\hbar
\left(U_\hbar(\g)'\right)_0+\left(U_\hbar(\g)'\right)_0^2\subset
\left(U_\hbar(\g)'\right)_0$ in the left hand side, and 
$\hbar^2 U_\hbar(\g) \subset \hbar U_\hbar(\g)$ in the right hand side,
\item[(c)] is the map
\begin{align*}
\hbar U_\hbar(\g)/\hbar^2 U_\hbar(\g)&\to
U_\hbar(\g)/\hbar U_\hbar(\g)\simeq U(\g)\\
x&\mapsto \hbar^{-1} x,
\end{align*}
\item[(d)] is induced by the pairing $\g^* \otimes \left(U(\g^*)^*\right)_0/
\left(U(\g^*)^*\right)_0^2 \to \KM$, quotient of the 
pairing
$\g^* \otimes \left(U(\g^*)^*\right)_0 \to \KM$, 
given by $\xi \otimes T \mapsto
\langle T,\canget(\xi)\rangle$ ($\canget$~: $\g^* \hookrightarrow U(\g^*)$ is
the canonical injection).
\end{description}
\end{proposition}
\dem
Let $\varphi'$ be an element of $ \left(U_\hbar(\g)'\right)_0$.
Recall that
$[\varphi']$ denotes its class
in $\left(U_\hbar(\g)'\right)_0/\big(\hbar \left(U_\hbar(\g)'\right)_0+
\left(U_\hbar(\g)'\right)_0^2\big)$.
Thanks to Lemma \ref{lemma:reduc}, we have 
$\left.{{\varphi'}\over{\hbar}}\right|_{U_{\hbar}\to U}\in
\g$ so $\hbox{(c)} \circ \hbox{(b)} \left([\varphi']\right)=
\cange\left(\left.{{\varphi'}\over{\hbar}}\right|_{U_{\hbar}\to U}\right)$.
Thus we should show that
\begin{equation}
\hbox{(d)}\circ\hbox{(a)} \left([\varphi']\right)=
\left.{{\varphi'}\over{\hbar}}\right|_{U_{\hbar}\to U}.
\label{eq:main}
\end{equation}
Let $S$ be a supplementary of $\g^*$ in 
$U(\g^*)_0$ (e.g., the image of $\oplus_{i \geq 2} S^i(\g^*)$ under
the symmetrization map). For $x \in \g$, let
$f_x \in U(\g^*)^*_0$ be defined by
$f_x(\xi)=\langle \xi,x \rangle$ for
$\xi \in \g^*$ and $f_x|_S=0$.
Then the class $[f_x]_0$ of $f_x$ in 
$U(\g^*)_0^*/\left(U(\g^*)_0^*\right)^2$ is independent of $S$.
Moreover, it is clear that we have the identity 
(d)$(f_x)=x$ for any $x \in \g$, so
$$\hbox{(d)}\left(\left[
f_{\left.{{\varphi'}\over{\hbar}}\right|_{U_{\hbar}\to U}}\right]_0\right)
=\left.{{\varphi'}\over{\hbar}}\right|_{U_{\hbar}\to U}.$$
So the identity (\ref{eq:main}) (and thus the proposition)
will be true if
\begin{equation}
\hbox{(a)}([\varphi'])=\left[f_{\left.{{\varphi'}\over{\hbar}}\right|_{U_{\hbar}\to
U}}\right]_0.
\label{eq:lemma}
\end{equation}
Both sides belong to
$U(\g^*)_0^*/\left(U(\g^*)_0^*\right)^2=\g$, so it suffices to show that their
pairings with any element
$\xi \in \g^*$ coincide. We have
\begin{align*}
\left\langle
\xi,\left[f_{\left.{{\varphi'}\over{\hbar}}\right|_{U_{\hbar}\to U}}\right]_0
\right\rangle&=\left\langle
\xi,\left.{{\varphi'}\over{\hbar}}\right|_{U_{\hbar}\to U}
\right\rangle, \\
\intertext{and}
\left\langle \xi, \hbox{(a)} ( [\varphi'])\right\rangle&=
\left\langle \xi, \left.{{\varphi'}\over{\hbar}}\right|_{U_{\hbar}\to U}\right\rangle
\hbox{ by construction of (a)},
\end{align*}
which proves identity (\ref{eq:lemma}).
\findem
\vskip1.1truecm

\centerline { ACKNOWLEDGEMENTS }

\vskip18pt

We would like to thank P. Etingof and P. Xu for discussions.

\end{document}